\documentclass[twoside,a4paper,10pt]{article}
\usepackage{amsmath,amsthm,amssymb,latexsym}
\usepackage[latin1]{inputenc}   

\usepackage[dvips]{graphicx}        

\newtheorem{thm}{Theorem}[section]
\newtheorem{theorem}[thm]{Theorem}
\newtheorem{corollary}[thm]{Corollary}
\newtheorem{lemma}[thm]{Lemma}

\theoremstyle{remark}

\begin{document}
\title{Hexagonal Tilings: Tutte Uniqueness}

\author{D. Garijo \footnote{ Dep. Matem\'{a}tica Aplicada I. Universidad de
Sevilla. Avda. Reina Mercedes s/n. 41012 Sevilla (Spain).
\{dgarijo,almar,pastora\}@us.es}, A. M\'{a}rquez $^*$, M.P.
Revuelta $^*$ }

\date{}

\maketitle

\thispagestyle{empty}

\begin{abstract}
We develop the necessary machinery in order to prove that hexagonal tilings are
uniquely
determined by their Tutte polynomial, showing as an example how to apply this technique to the
toroidal hexagonal tiling.
\end{abstract}

\section{Introduction}

The main result in this paper has to do with graphs determined by their
Tutte polynomial. This is
a two variable polynomial
$T(G;x,y)$ associated with any graph $G$,
which contains interesting information about $G$. For instance,
$T(G;2,0)$ is the number of acyclic orientations of $G$ and $T(G; x,0)$
with $x<0$ is the
chromatic polynomial associated with $G$. As a natural extension
of the concept of chromatically
unique graph, the notion of Tutte unique graph was studied in \cite{miernoy}.
A graph $G$ is said to
be \emph{Tutte unique} if $T(G;x,y)=T(H;x,y)$ implies $G\cong H$ for any
other graph $H$. A
common topic in the study of this invariant is the search of large families
of Tutte
unique graphs. In 2003, Garijo, Márquez, Mier, Noy, Revuelta \cite{gmr},\cite{mmnr}
 present locally grid
graphs as the first large family of graphs uniquely determined by their
Tutte polynomial. In this
paper we prove that the toroidal hexagonal tilings also
satisfy this uniqueness property.
The technique developed here can also be applied to the rest
of hexagonal tilings and their dual
graphs, the locally $C_6$ graphs. The main motivation for this study
is the search of a
generalization to locally plane graphs of the well known
relationship existing between the Tutte
polynomials of plane graphs and their duals \cite{tutte}.

\begin{figure}[htb]
\begin{center}
\includegraphics[width=10mm]{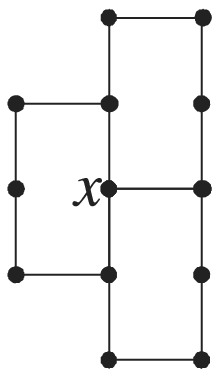} 
\end{center}\caption{Hexagonal structure around $x$}
\label{F0}
\end{figure}

A hexagonal tiling, $H$, is defined as a connected cubic graph of girth 6,
having a collection of
$6$-cycles, $\textit{C}$, such that every $2-$path is contained
in precisely one cycle of
$\textit{C}$ ($2-$path condition). A hexagonal tiling is simple
(that is, without loops and
multiple edges) and every vertex belongs to exactly three hexagons
(cycles of length 6)(Figure
\ref{F0}). Every hexagon of the tiling is called a \emph{cell}.

Let $G=(V,E)$ be a graph with vertex set $V$ and edge set $E$.
The \emph{rank} of a subset
$A\subseteq E$ is $r(A)=|A|-k(A)$, where $k(A)$
is the number of connected components of the
spanning subgraph $(V,A)$. The Tutte polynomial is defined as
follows:
$$T(G;x,y)=\sum_{ A\subseteq E} (x-1)^{r(E)-r(A)}(y-1)^{|A|-r(A)}$$

$T(G;x,y)$ contains exactly the same information about $G$
as the rank-size generating polynomial which is defined as:
$$R(G;x,y)=\sum_{ A\subseteq E}x^{r(A)}y^{|A|}$$

\noindent The coefficient of $x^iy^j$ in $R(G;x,y)$ is the number of spanning
subgraphs in $G$ with rank $i$ and $j$ edges, hence the
Tutte polynomial tell us for every $i$ and
$j$ the number of edges-sets in $G$ with rank $i$ and size $j$.

The main strategy in the paper is the following.
We first recall some definitions and results from
\cite{ggmr, thomassen}, such as, the classification theorem of hexagonal tilings.
In Section 3, we
show that these graphs are locally orientable,
using the relationship existing between hexagonal
tilings and locally grid graphs \cite{ggmr}. The local orientation
of these graphs allows us to
count edge-sets and it is one of the basic tools used to prove the Tutte uniqueness results.
Thus, as an example,
in Section 5 we apply the technique developed to show that for
every hexagonal tiling $H$ not
isomorphic to $H_{k,m,0}$ with $2k(m+1)$ vertices there is at
least one coefficient of the
rank-size generating polynomial in which both graphs differ.

\

\section{Preliminary Results}

In this section we state the classification theorem of hexagonal tilings proved in \cite{ggmr}.
There exits an extensive literature on this topic. See for instance the works done by Altshuler
\cite{altshuler1, altshuler2}, Fisk \cite{fisk1, fisk2} and Negami \cite{negami1,negami2}. In
\cite{ggmr}, following up the line of research given by Thomassen \cite{thomassen}, we add two new
families to the classification theorem given in \cite{thomassen} proving that with these families
we exhaust all the cases.

\begin{theorem}\label{clasificacion}
\cite{ggmr} If $G$ is a hexagonal tiling with $N$ vertices, then one and only one of the following
holds: $$\begin{tabular}{lcl}
 {\bf A)} & $G\simeq H_{k,m,r}$ & with  $N=2k(m+1)$, $0\leq r\leq \lfloor k/2 \rfloor $, $m\geq
2$, $k\geq 3$. If $m=1$ then $k>3$  \\ \noalign{\smallskip} & & and $\lfloor k/2 \rfloor \geq
r\geq 2 $. If $m=0$ then $k>3$ and $3\leq r\leq k$  \\ \noalign{\medskip} {\bf B)} & $G\simeq
H_{k,m,a}$ &
with  $N=2k(m+1)$, $m\geq 2$, $k\geq 3$. \\
\noalign{\medskip} {\bf C)} & $G\simeq H_{k,m,b}$ & with
$N=2k(m+1)$, $k$ even, $m$ odd, $m\geq 3$, $k\geq 4$. \\
\noalign{\medskip}
 {\bf D)} & $G\simeq H_{k,m,c}$ & with  $N=2k(m+1)$, $m\geq 1$, $k$ even, $k\geq 6$. \\ \noalign{\medskip}
{\bf E)} & $G\simeq H_{k,m,f}$ & with  $N=2k(m+2)$, $k$ odd, $m\geq 0$, $k\geq 7$. \\
\noalign{\medskip} {\bf F)} & $G\simeq H_{k,m,g}$ &
with  $N=2(m+1)(k+2)$, $k\geq m+1$, $m\geq 3$. \\
\noalign{\medskip} {\bf G)} & $G\simeq H_{k,m,h}$ & with $N=2(m+1)(k+1)$, $k< m-1$, $k\geq 2$.
\\ \end{tabular}$$
\end{theorem}

\begin{figure}[htb]
\begin{center}
\includegraphics[width=130mm]{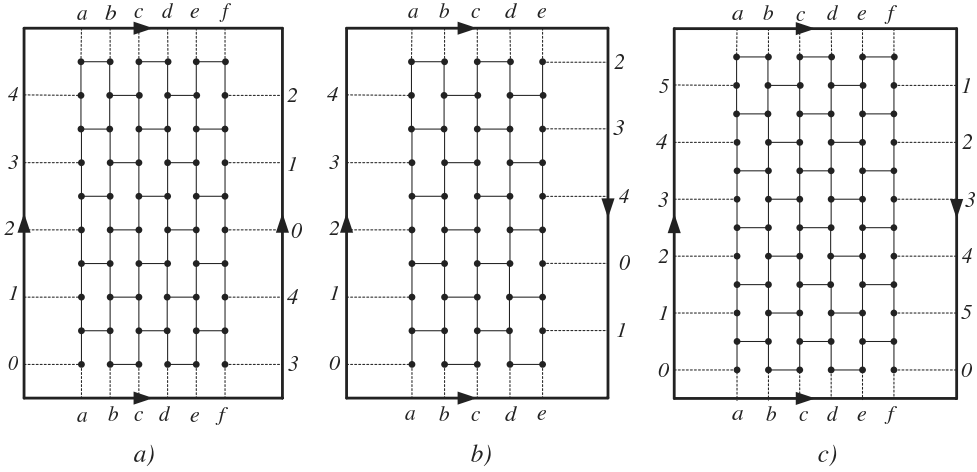} 
\end{center}
\caption{a) $H_{5,5,r}$ b) $H_{5,4,a}$ c) $H_{6,5,b}$}
\label{Hrab}
\end{figure}

\begin{figure}[htb]
\begin{center}
\includegraphics[width=100mm]{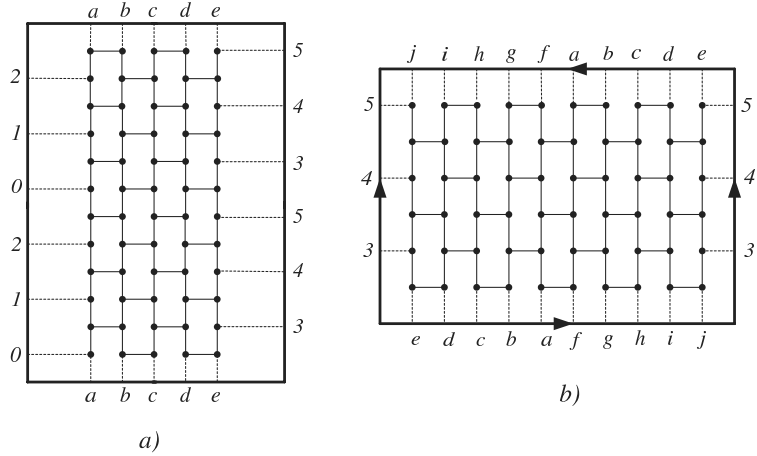} 
\end{center}
\caption{a) $H_{6,4,c}$  b) Embedding of $H_{6,4,c}$ in the Klein
bottle } \label{Hc}
\end{figure}

\begin{figure}[htb]
\begin{center}
\includegraphics[width=90mm]{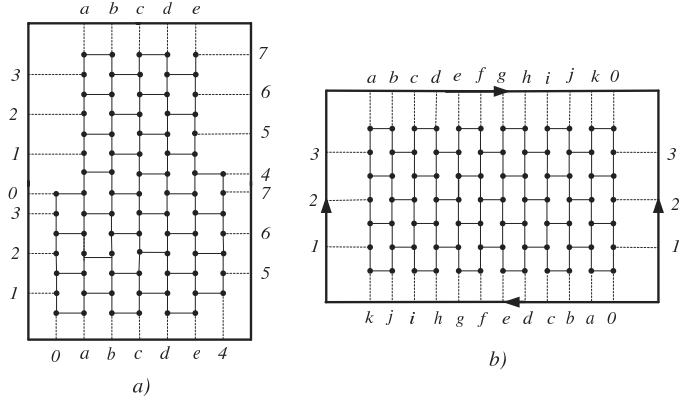} 
\end{center}
\caption{a) $H_{7,4,f}$  b) Embedding of $H_{7,4,f}$ in the Klein
bottle} \label{Hf}
\end{figure}

\begin{figure}[htb]
\begin{center}
\includegraphics[width=110mm]{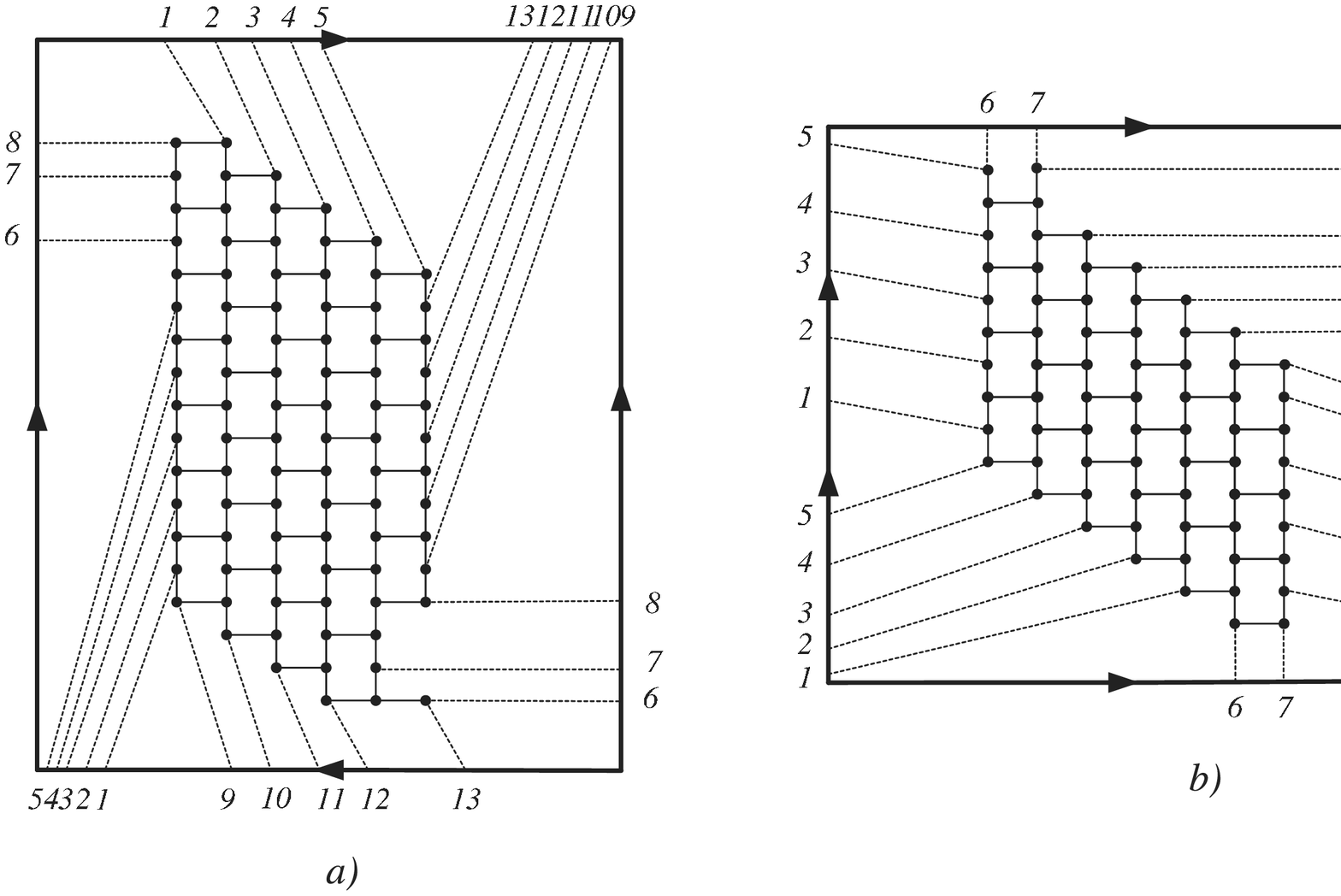} 
\end{center}
\caption{a) $H_{7,4,g}$  b)  $H_{6,4,h}$} \label{Hgh}
\end{figure}

Some examples of hexagonal tilings are given in Figures \ref{Hrab}, \ref{Hc}, \ref{Hf} and
\ref{Hgh}. The parameters $k$ and $m$ fix, respectively, the height and the breath of the
structures, called \emph{hexagonal wall of length $k$ and breath $m$} (Figures \ref{Hrab},
\ref{Hc}, \ref{Hf}) and \emph{hexagonal ladder of length $k$ and breath $m$} (Figure \ref{Hgh}),
in which it is added edges to construct the different families of hexagonal tilings (see
\cite{ggmr}).

Given two cycles $C$ and $C'$ in a hexagonal tiling $G$, we say that $C$ is \emph{locally
homotopic} to $C'$ if there exists a cell, $H$, with $C\cap H $ connected and $C'$ is obtained
from $C$ by replacing $C\cap H $ with $H-(C\cap H)$. A \emph{homotopy} is a sequence of local
homotopies. A cycle in $G$ is called \emph{essential} if it is not homotopic to a cell. Otherwise
it is called \emph{contractible}. This definition is equivalent to the one given for a graph
embedded in a surface \cite{tesis}. Let $l_{G}$ be the minimum length of the essential cycles of
$G$, $l_{G}$ is invariant under isomorphism. An edge-set contained in a hexagonal tiling is called
\emph{a normal edge-set} if it does not contain any essential cycle.

\begin{figure}[htb]
\begin{center}
\includegraphics[width=60mm]{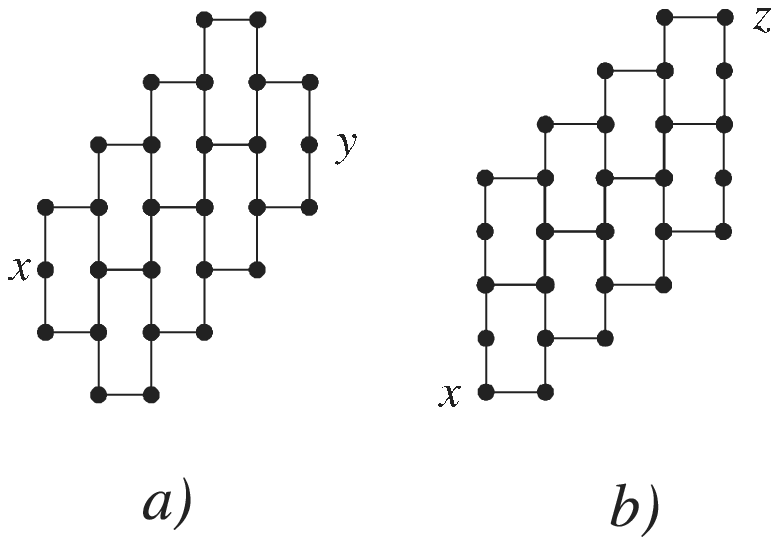} 
\end{center}\caption{ a) Displaced ladder $4\times 2$  b) Ladder $4\times
2$}\label{F12}
\end{figure}

There are two different ways of pasting together $j$
ladders each one containing $i$ hexagons,
from which we obtain two structures, called the \emph{ladder $i\times j$} and the \emph{displaced
ladder $i\times j$ }, shown in Figure \ref{F12}. Note that
the number of shortest paths between
$x$ and $y$ or between $x$ and $z$ in a
ladder $i\times j$ or in a displaced ladder $i\times j$
is $\left(\begin{array}{c} i+j \\
j\end{array}\right)$ and the length of these paths is $2(i+j)-1$.

Every hexagonal tiling is obtained by adding edges to a hexagonal wall or to a hexagonal ladder
(except for $H_{k,m,h}$, in which we have also added two vertices) \cite{ggmr}. These edges are
called \emph{exterior edges} and every essential cycle must contains at least one of these edges (
the edges $\{(0,2k+m), (0,2k+m-1)\}$ and  $\{(m-k-1,3k+2), (m-k-1,3k+1)\}$ in $H_{k,m,h}$ are not
considered exterior edges).

\

\section{Local Orientability of Hexagonal Tilings}

In this section we prove that hexagonal tilings are locally orientable.
This fact allow us to
count edge-sets and it is the tool to distinguish at least one coefficient
of the Tutte
polynomials associated to two non isomorphic hexagonal tilings.
We first recall a minor
relationship existing between hexagonal tilings and locally grid graphs (see \cite{ggmr}).

\begin{figure}[htb]
\begin{center}
\includegraphics[width=10mm]{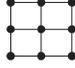} 
\end{center}
\caption{Locally Grid Structure}\label{grid0}
\end{figure}

We say that a $4-$regular, connected graph $G$ is a \emph{locally
grid graph} if for every vertex $x$ there exists an ordering
$x_1,x_2,x_3,x_4$ of $N(x)$ and four different vertices
$y_1,y_2,y_3,y_4$, such that, taking the indices modulo 4,
$$\begin{tabular}{ccl}
$N(x_i)\cap N(x_{i+1}) $ & $=$ & $\{x,y_i\}$ \\
\noalign{\medskip} $N(x_i)\cap  N(x_{i+2}) $ & $=$ & $\{x\}$ \\
\end{tabular}$$ \noindent and there are no more adjacencies among $\{
x, x_1, \ldots ,x_4,y_1, \ldots , y_4\}$ than those required by
this condition (Figure \ref{grid0}).

Every locally grid graph is a minor of a hexagonal tiling and in \cite{ggmr}
it was proved that
there exits a biyective minor relationship preserved by duality between
hexagonal tilings with the
same chromatic number and locally grid graphs.
A perfect matching $P$ is selected in each one
of the families given in the classification theorem of hexagonal tilings,
except in one case,
$H_{k,m,f}$ in which the set of selected edges is not a
 matching. $P$ contains two edges of each hexagon. In $H_{k,m,f}$ there are $k-1$
hexagons in which we select four edges pairwise incident. A locally grid graph is obtained by
contracting the edges of this matching, and deleting parallel edges if necessary. There are just
two cases in which parallel edges must be deleted, $H_{k,m,f}$ and $H_{k,m,h}$. These perfect
matchings and the set of edges of $H_{k,m,f}$  verify that if we have two hexagonal tilings with
the same chromatic number, the result of the contraction of their perfect matchings (or the
selected edges in $H_{k,m,f}$) are two locally grid graphs belonging to different families. This
condition is going to be essential in the search of a generalization, for these families, of the
relationship existing between the Tutte polynomial of a plane graph and its dual \cite{tutte}.

In \cite{mmnr} it was proved that locally grid graphs are locally orientable. Given a vertex $v$
of a locally grid graph $G$, two edges incident with $v$ are said to be \emph{adjacent} if there
is a square containing them both; otherwise they are called \emph{opposite}. An orientation at a
vertex $v$ consists of labeling the four edges incident with $v$ bijectively with the labels $N$,
$S$, $E$, $W$ in such a way that the edges labelled $N$ and $S$ are opposite, and so are the ones
labelled as $E$ and $W$.

\begin{lemma}\label{localorienta}
Hexagonal tilings are locally orientable.
\end{lemma}

\begin{proof}
Let $H$ be a hexagonal tiling. By Theorem \ref{clasificacion},
$H$ is isomorphic to one of the following graphs: $H_{k,m,r}$,
$H_{k,m,a}$, $H_{k,m,b}$, $H_{k,m,c}$, $H_{k,m,h}$, $H_{k,m,g}$,
$H_{k,m,f}$.

Suppose first that $H$ is not isomorphic to $H_{k,m,f}$. For
every vertex $v\in V(H)$ there is exactly one edge of the perfect
matching $P$ incident with $v$, call it $\{v,v'\}$. Consider the
union of the local structures of $v$ and $v'$, labeling the other
two edges incident with $v$ as $e_1$ and $e_2$ and the other two
edges incident with $v'$ as $e_3$ and $e_4$. Contract all the
edges of $P$ belonging to this union, deleting parallel edges if
necessary. As it is shown in Figure \ref{localorientation}a we
obtain a grid $3\times 3$, hence we can label the edges $e_1$,
$e_2$, $e_3$ and $e_4$ bijectively with the labels $N$, $S$, $E$,
$W$ following the criterion established for the locally grid
graphs. The edge $\{v,v'\}$ is considered opposite to the unique
edge incident with $v$ and belonging to the perfect matching,
$P'$, formed by the edges $\{(i,l) (i,l+1)\}$, $0\leq i\leq m$ not
belonging to $P$.

\begin{figure}[htb]
\begin{center}
\includegraphics[width=70mm]{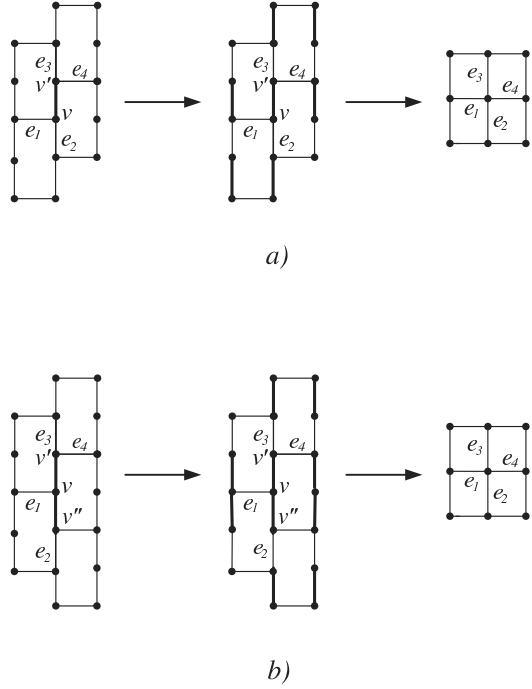} 
\end{center}
\caption{a) Local orientation at $v\in V(H)$ not isomorphic to
$H_{k,m,f}$ b) Local orientation at $v\in V(H_{k,m,f})$}
\label{localorientation}
\end{figure}

If $H$ is isomorphic to $H_{k,m,f}$, then there are either one or two edges belonging to $P$ and
incident with $v$. If $v$ does not belong to a hexagon that has more than two edges of $P$, we
follow the same process as in the previous case. Suppose, now, that $v$ is a vertex belonging to
one of the hexagons that contains four edges belonging to $P$. We consider two cases:

{\bf (1)} There are two edges of $P$ incident with $v$, call them
$\{v,v'\}$ and $\{v,v''\}$. The union of the hexagonal structures
of $v$, $v'$ and $v''$ is composed by five hexagons. Label the
unique edge incident with $v$ and not belonging to $P$ as $e_1$.
The edge incident with $v''$ not belonging to $P$ and sharing an
hexagon with $e_1$ is labelled $e_2$. Finally, the edges incident
with $v'$ no belonging to $P$ are labelled $e_3$ and $e_4$.
Contracting all the edges of $P$ belonging to this union and
deleting parallel edges, we also obtain a grid $3\times 3$ and we
can label the edges incident with $v$, $v'$ and $v''$ bijectively
with the labels $N$, $S$, $E$, $W$ (Figure
\ref{localorientation}b). $\{v,v'\}$ and $\{v,v''\}$ are
considered opposite to the unique edge incident with $v'$ and
$v''$ respectively, belonging to $P'$.

{\bf (2)} There is one edge of $P$ incident with $v$, call it
$\{v,v'\}$. There has to be one edge, different from $\{v,v'\}$,
incident with $v'$ and belonging to $P$, call it $\{v',v''\}$.
The union of the hexagonal structures of $v$ and $v''$ is
composed by five hexagons. The two edges incident with $v$ not
belonging to $P$ are labelled $e_1$ and $e_2$. The unique edge
incident with $v'$ not belonging to $P$ is labelled $e_3$, and the
edge incident with $v''$ not belonging to $P$ or to the hexagon
that contains four edges of $P$, is labelled $e_4$. Contracting
all the edges of $P$ belonging to this union and deleting
parallel edges, we also obtain a grid $3\times 3$ and we can
label the edges incident with $v$, $v'$ and $v''$ bijectively
with the labels $N$, $S$, $E$, $W$. $\{v,v'\}$ and $\{v',v''\}$
are considered opposite to the edges incident with $v$ and $v''$
respectively, belonging to the hexagon that contains four edges
of $P$.

\end{proof}

The previous proof give us a definition of \emph{adjacent and
opposite edges} for hexagonal tilings. Two edges incident with a
vertex $v$ and no belonging to $P$ are called \emph{adjacent} or
\emph{opposite} if they are so in the locally grid structure
resulting from applying the operations explained above.
Two edges incident with $v$ and belonging to $P$ are
always opposite.

By the orientation $(v,v',e,f)$ we mean that $v$ is the origin
vertex, $\{v,v'\}\in P$, $e$ is the unique edge incident with $v$
and opposite to the other two edges incident with $v$. This edge
never belongs to $P$ and it is labelled $E$. Finally, $f$ is an
edge incident with $v$ and labelled $N$. The
orientation $(v,v',v'',e,f)$ means that $v$ is the origin vertex,
$\{v,v'\}, \{v, v''\}\in P$, $e$ is the unique edge incident with
$v$ and opposite to the other two edges incident with $v$, and
$f$ is an edge incident with $v$ labelled $N$. We will denote
both cases as $(v,v',e,f)$.

If $\alpha\in \{N,S,E,W\}$, we denote $\alpha^{-1}$ the label
opposite to $\alpha$. If we fix an orientation at $v\in V(H)$
($v$ is the origin vertex) then every vertex $u$ adjacent to $v$
is unambiguously oriented, because the orientation at $v$ induces
an orientation at $u$: if $\{u,v\}$ is labelled $\alpha$ from
$v$, it is labelled $\alpha^{-1}$ from $u$. If $tzyxvw$ is a
hexagon and $\{x,v\}$, $\{w,v\}$ are labelled $\beta$ and $\alpha$
from $v$ respectively, then $\{x,y\}$ is $\beta$ from $x$,
$\{z,y\}$ is $\alpha$ from $y$, $\{z,t\}$ is $\beta^{-1}$ from
$z$ and $\{t,w\}$ is $\beta^{-1}$ from $t$. We also can translate
the orientation at $v$ to all the vertices of a path beginning at
$v$, but two different paths joining $v$ and $u$ can induce
different orientation at $u$. This does not happen if the union
of these two paths is a contractible cycle. In fact, we can
transform one path into the other one using three elementary
transitions and their inverse: if $e$, $f$, $g$, $h$, $i$, $j$
are the edges of a hexagon ordered cyclically, we can change $e$,
$f$, $g$, $h$, $i$ by $j$; $e$, $f$, $g$, $h$ by $i$, $j$ and
$e$, $f$, $g$ by $h$, $i$, $j$ and these operations do not change
the orientation at the endpoint.

Let $A$ be a connected normal edge set in a hexagonal tiling $H$.
Fixed an orientation at a vertex $v\in V(A)$, then all the
vertices in $V(A)$ are unambiguously oriented. Every path in $A$
can be described as a sequence of the labels $\{N,S,E,W\}$
(Figure \ref{orientaconjuntos}a). This allow us to assign
coordinates to every vertex in $V(A)$: $v$ has coordinates
$(0,0)$ and $u\in V(A)$ has coordinates $(i,j)$ if in one path in
$A$ joining $v$ to $u$ (taking into account that the orientation
at the endpoint induced by different paths joining $v$ to $u$
does not change), $i$ is the number of labels $E$ minus the
number of labels $W$ and $j$ is the number of labels $N$ minus
the number of labels $S$.

\begin{lemma}\label{asignacion}
Let $A$ be a connected normal edge set in a hexagonal tiling $H$
with $|A|\leq l_H+3$, then different vertices in $V(A)$ have
different coordinates.
\end{lemma}

\begin{proof}
First, we are going to prove that with a fixed orientation
$(v,v',e,f)$ no vertex except the origin can have coordinates
$(0,0)$. Suppose that there exists $x\in V(A)$ such that $x$ has
coordinates $(0,0)$, then there exits a path $R$ joining $v$ to $x$
with as many labels $N$ as labels $S$ and as many labels $E$ as labels $W$. We
are going to show that $x=v$ by induction on the length of this
path.

If $|R|=6$ and $v$ is the origin vertex, the first labels of $R$
can be: $NN$, $NW$, $E$, $SS$, or $SW$ and using the elementary
transitions described previously, we can change these labels by
$ENNW$, $SWNN$, $NNESS$, $ESSW$ and $NWSS$ respectively. To
obtain as many labels $N$ as labels $S$ and as many labels $E$ as labels $W$ it
must be that $v=x$. We assume that every path $R$ (not a cycle) with
$|R|<n$ does not have as many labels $N$ as labels $S$ and as many
labels $E$ as labels $W$ and we prove the result for $|R|=n$. Suppose
that there exits a path $R$ joining $v$ and $x$, not a cycle, with
$|R|=n$ and having as many labels $N$ as those labelled $S$ and as many
$E$ as $W$, then
there exits a subsequence labelled $\beta (\alpha^{-1})^l \beta
\alpha^s \beta^{-1} (\alpha^{-1})^m  \beta^{-1}$ or $\beta
\alpha^l \beta \alpha^s \beta^{-1} \alpha^m \beta^{-1}$ (or other
analogous subsequence). In the first case we can change this
subsequence by $\alpha^{s-l-m}$ keeping the orientation at the
endpoint. Analogously, in the second case the subsequence can be
changed by $\alpha^{s+l+m}$. In any case, the change gives rise to
a new path $R'$ joining $v$ and $x$ with $|R'|<n$, in a new
connected edge set $A'$. This set is also a normal edge set
because $|A'|\leq |A|-4\leq l_H-1$. By hypothesis, $R'$ does not
have as many labels $N$ as $S$ and $E$ as $W$, hence neither does
$R$ and we obtain a contradiction.

Given two different vertices, $u$ and $w$ in $V(A)$, we consider
the paths $R_u$ and $R_w$ joining $u$ to $v$ and $w$ to $v$
respectively. Since $A$ is connected, there is a path $R_w \cup R$
joining $v$ to $u$ and keeping the orientation at $u$. Hence,
it is easy to prove that $u$ and $w$ can not have the same
coordinates.
\end{proof}

From now on, and unless otherwise stated, we consider that all normal edge-sets have at most
$l_H+3$ elements.

\section{Counting Edge-Sets}

The aim of this section is to find a system which allow us to
count edge-sets in order to prove that there is at least one
coefficient of the rank-size generating polynomial in which two
non isomorphic hexagonal tilings differ. Our first step is to
codify the edge-sets of a hexagonal tiling.

Let $L^{\infty}$ be the infinite plane square lattice, that is,
the infinite graph having as vertices $\mathbb{Z} \times
\mathbb{Z}$ and in which $(i,j)$ is joined to $(i-1,j), (i+1,j),
(i,j-1), (i,j+1)$. We delete the edges $\{(i,2j), (i-1,2j)\}$ if
$i$ even and $\{(i,2j+1), (i+1,2j+1)\}$ if $i$ odd, obtaining a
\emph{infinite plane hexagonal tiling $H^{\infty}$}. Let $\Omega$
be the group of graph automorphisms of $H^{\infty}$. This
group is generated by translations, symmetries and rotations
of the plane that map vertices to vertices. Let
$\Sigma(H^{\infty})$ be the set of all finite connected edge-sets
of $H^{\infty}$. An equivalence relation is defined in
$\Sigma(H^{\infty})$ as follows:

$$B_1\sim B_2 \Leftrightarrow
\exists \sigma \in \Omega, \sigma(B_1)=B_2$$

If  $\Lambda$ is a set of representatives of $\Sigma(H^{\infty})/\sim$ such that
 every $B\in \Lambda$ contains the vertex $(0,0)$, then $\Lambda$
 covers all the possible shapes that a normal edge-set could have.

 We are going to define a set $\Gamma$ of words over the alphabet
 $\{N,S,E,W\}$ that represents all edge-
 sets of $\Sigma(H^{\infty})/\sim$ . Label $\{(0,0),(0,1)\}$ as $N$,
 $\{(0,0),(1,0)\}$ as $E$, $\{(0,1),(-1,1)\}$ from $(0,1)$ as $W$ and
 $\{(0,0),(-1,0)\}$ as $S$. For every $B\in \Lambda$ there is a
 sequence $\gamma_B$ over the the alphabet
 $\{N,S,E,W\}$ such that beginning at $(0,0)$ and following the
 instructions given by $\gamma_B$ the unique edges covered are
 those of $B$. Note that one edge can be covered more than once
 and that the sequence $\gamma_B$ is not unique (Figure \ref{orientaconjuntos}a).
Since every sequence over the alphabet $\{N,S,E,W\}$ can not represent a
 set of $\Lambda$ (for instance, a sequence with two consecutive labels $E$)
we are going to establish some restrictions on this
 alphabet:

 {\bf (1)} No sequence begins with the label $W$.

 {\bf (2)} There are no two consecutive labels $E$ (analogously $W$).

 {\bf (3)} If a sequence begins with a label $N$ (analogously $S$),
there must be an odd number of
 labels $N$ and $S$ before having a label $W$; and an even number
 if the label is $E$. Hence, there are an odd number of
 labels $N$ and $S$ between two labels $W$ and an even number
 between two labels $E$.

{\bf (4)} If a sequence begins with a label $E$, there are an odd number of
 labels $N$ and $S$ between two labels $E$ and an even number between two labels $W$.

 Call $\Gamma$ the set of words over the alphabet $\{N,S,E,W\}$
 verifying the three previous conditions. These restrictions do
 not modify the fact that every $B\in \Lambda$ has a word
 $\gamma_B$ associated. Conversely, given $\gamma\in \Gamma$ there exits $B\in
 \Lambda$, such that $\gamma= \gamma_B$. It will be denoted
 $B(\gamma)$.

 The next step is to assign one word from $\Gamma$ to each normal
 connected edge-set of a hexagonal tiling, $H$. Fix an orientation
$(v,v',e,f)$ at $H$, being
 $v\in V(H)$ the origin vertex. Following the code given by
 $\gamma$ from $v$ with orientation $(v,v',e,f)$, an edge
 set $A^H(\gamma, v,v',e,f)$ is obtained, called the \emph{instance of $\gamma$ in
 $H$} and we call $(v,v',e,f)$ an \emph{orientation of the
 instance}. An instance of $\gamma$ can have different
 orientations associated to it, the number of them depends only on the
 symmetries of $B(\gamma)$ and it is called $sym(\gamma)$. The \emph{length of
 $\gamma$} is the number of labels $N$, $S$, $E$, $W$ that give
 rise to $\gamma$.

 \begin{lemma}\label{Bgamma}
Let $H$ be a hexagonal tiling and $\gamma\in \Gamma$. If
$r(B(\gamma))=l_H-1$ and $|B(\gamma)|\leq l_H$, then $A^H(\gamma,
v,v',e,f)$ is a normal edge-set and it has the same rank and size
as $B(\gamma)$.
 \end{lemma}

\begin{proof}
We are going to prove this result by induction on the length of
$\gamma$. If $\gamma$ has length equal to one, the result is
trivial. Let $\gamma'$ be the word resulting from removing the last
label from $\gamma$. By hypothesis $A^H(\gamma',v,v',e,f)$ is
normal and it has the same rank and size as $B(\gamma')$. When
we add the last label of $\gamma$ to $\gamma'$, there are three
possibilities on $B(\gamma')$:

(1) We are adding one edge between two vertices belonging to
$B(\gamma')$

(2) We are adding an isthmus.

(3) We do not add any edge.

\

{\bf (1)} Let $x'$ be the last vertex of $A^H(\gamma',v,v',e,f)$
covered by $\gamma'$. Since $A^H(\gamma',v,v',e,f)$ is normal all
the vertices of $A^H(\gamma',v,v',e,f)$ are unambiguously
oriented. Let $c$ be the added edge. This edge is labelled
$\alpha$ from $x'$ following the orientation $(v,v',e,f)$.
Suppose that $A^H(\gamma,v,v',e,f)$ is not a normal edge-set,
then when we have added the edge $c$ we have created an essential
cycle. Since $|B(\gamma)|\leq l_H$, $A^H(\gamma,v,v',e,f)$ has to
be an essential cycle. $B(\gamma)$ contains at least one cycle
because it has been obtained joining two vertices of
$B(\gamma')$. Hence $\gamma$ contains a subword that can be
$\beta (\alpha^{-1})^s \beta^{-1}\alpha^s$, $\beta
(\alpha^{-1})^s \beta (\alpha^{-1})^t \beta^{-1} \alpha^r
\beta^{-1}$, $\beta (\alpha^{-1})^s \beta (\alpha^{-1})^t
\beta^{-1} (\alpha^{-1})^r \beta^{-1}$ (or other analogous cases
to the last ones) (Figure \ref{orientaconjuntos}b). We change
these subwords by $\beta (\alpha^{-1})^{s-2}
\beta^{-1}\alpha^{s-2}$, $(\alpha^{-1})^{t+s-r}$,
$(\alpha^{-1})^{t+s+r}$ respectively.

\begin{figure}[htb]
\begin{center}
\includegraphics[width=70mm]{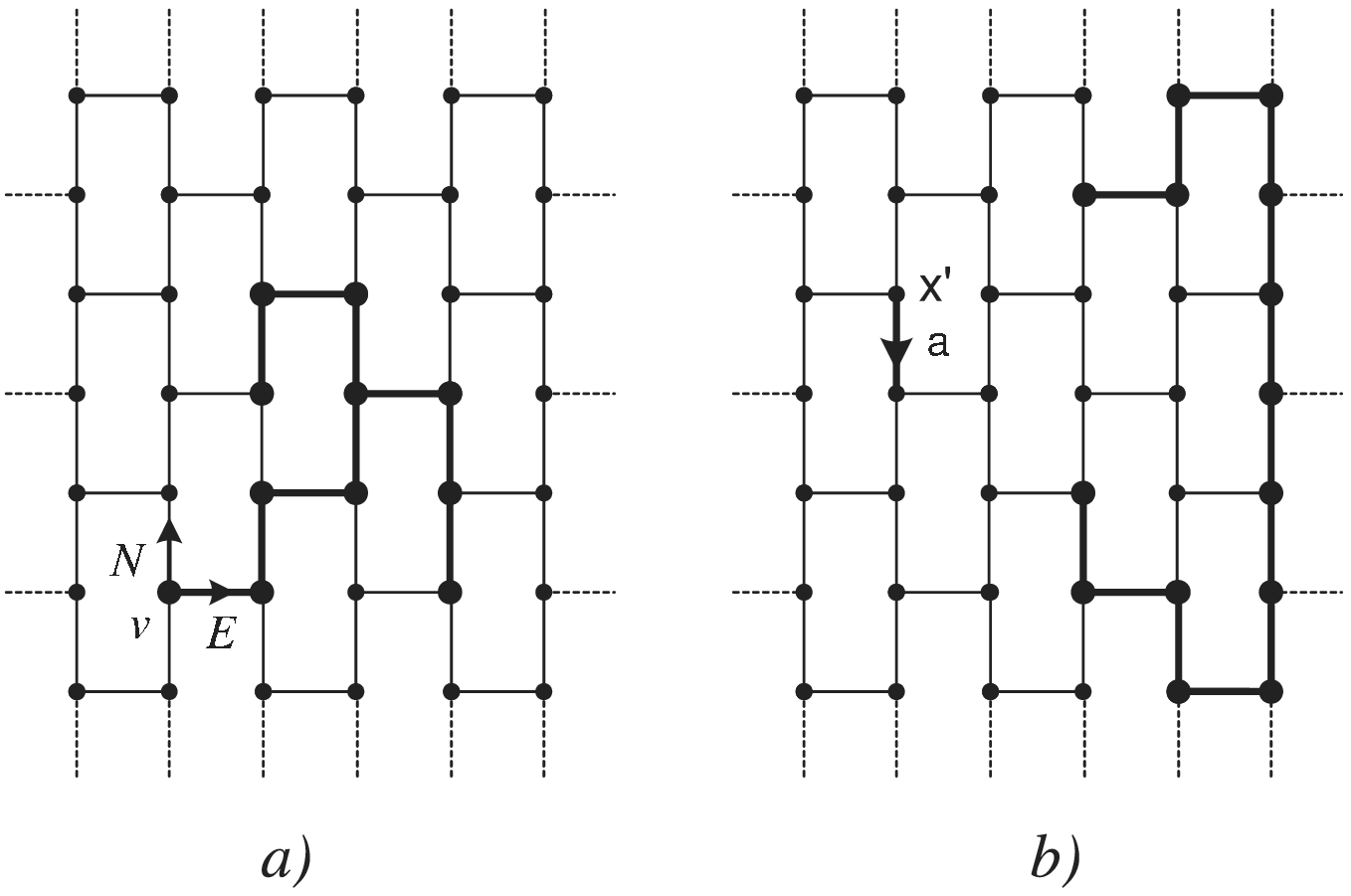} 
\end{center}
\caption{a) A path described using two different sequence over
$\{N,S,E,W\}$, the first one is $ENENNWSNESESS$ and the second
one is $ENENESSNNWNWS$ b) The marked path labelled as
$abab(a^{-1})^6 b^{-1} a b^{-1}$ is changed by $(a^{-1})^3$}
\label{orientaconjuntos}
\end{figure}

\noindent In any case, a word $\gamma^*$ is obtained of length
smaller than the length of $\gamma$, hence
$A^H(\gamma^*,v,v',e,f)$ is a normal edge-set and contains a
cycle. If we consider the simply connected region determined by
$A^H(\gamma^*,v,v',e,f)$ and we add, depending on the case, one
hexagon or two columns of hexagons, a simply
connected region is obtained which is determined by $\gamma$. Since
$A^H(\gamma,v,v',e,f)$ is an essential cycle, at least one vertex
of the added hexagons belongs to $A^H(\gamma^*,v,v',e,f)$.
Therefore, $A^H(\gamma,v,v',e,f)$ has at least two cycles and we
obtain a contradiction. The equality of ranks
and sizes remains to be proved. By hypothesis of induction
$r(A^H(\gamma',v,v',e,f))=r(B(\gamma'))$ and
$|A^H(\gamma',v,v',e,f)|=|B(\gamma')|$. Furthermore
$|B(\gamma)|=|B(\gamma')|+1$ and $r(B(\gamma))=r(B(\gamma'))$.
Since  $A^H(\gamma,v,v',e,f)$ is normal,
$r(A^H(\gamma,v,v',e,f))=r(A^H(\gamma',v,v',e,f))$ and
$|A^H(\gamma,v,v',e,f)|=|A^H(\gamma',v,v',e,f)|+1$.

\

{\bf (2) and (3)} The second case is easier than the previous
one. Using a similar reasoning we obtain that if
$A^H(\gamma,v,v',e,f)$ is not a normal edge-set, then it is an
essential cycle. But an essential cycle only could be created
joining two vertices belonging to $B(\gamma')$. The equality of
ranks and sizes is proved as in case (1). The third case is
trivial because if we do not add any edge,
$A^H(\gamma,v,v',e,f)=A^H(\gamma',v,v',e,f)$.
\end{proof}

\begin{lemma}\label{asignarBaA}
Let $A\subseteq E(H)$ be a connected normal edge-set. Fixed an
orientation $(v,v',e,f)$ with $v\in V(A)$, there exits a unique
edge-set $B\subseteq E(H^{\infty})$ and a unique graph
isomorphism $\varphi:V(A)\rightarrow V(B)$ such that:

(1) $\varphi(v)=(0,0)$.

(2) If the edge $\{x,y\}\in A$ is labelled $\alpha$ from $x$,
then $\{\varphi(x), \varphi(y)\} $ is labelled $\alpha$ from
$\varphi (x)$ according to the orientation $((0,0),
\{(0,0),(1,0)\}, \{(0,0),(0,1)\}$.
\end{lemma}

\begin{proof}
Fixed an orientation $(v,v',e,f)$ with $v\in V(A)$, we can assign
coordinates to the vertices of $A$. Since $A$ is normal, this
assignment is unambiguous. We define $\varphi:V(A)\rightarrow
V(H^{\infty})$ by $\varphi (x)=(i,j)$ if in every path in $A$
joining $v$ to $x$, $i$ is the number of labels $E$ minus the
number of labels $W$ and $j$ is the number of labels $N$ minus
the number of labels $S$. Furthermore, we define $\bar{\varphi}:
A\rightarrow E(H^{\infty})$ by $\bar{\varphi}
(\{x,y\})=\{\varphi(x) ,\varphi(y)\} $. Taking
$B=\bar{\varphi}(A)$, it is easy to prove that $\varphi$ is a
graph isomorphism and that $\bar{\varphi}$ is the isomorphism
induced on edges by $\varphi$. The statements (1) and (2) of the
lemma are consequence of the assignment of coordinates following
an orientation $(v,v',e,f)$ in a hexagonal tiling $H$, and the
assignment of labels in the edges incident with the origin vertex
$v$. We prove the uniqueness of $B$ by induction on
$m=$max$\{d(v,x), x\in V(A)\}$. Suppose that there exits a pair
$(\varphi', B')$ verifying the statement of the lemma. If $m=1$,
$d(v,x)=1$ $\forall x\in V(A)$ with $x\neq v$ hence $|A|\leq 3$.
By (1) and (2), it is immediate to show that $B=B'$. Assume the
result for $m\leq l$ and we prove it for $m=l+1$. There exits
$x_1, \ldots, x_n \in V(A)$ such that $d(v,x_i)=l+1$ and
$d(v,x)<l+1$ for all $x\neq x_i$. By induction on $n$ we show
that $B$ is unique for $m=l+1$. If $n=1$, let $x_1$ be the unique
vertex such that $d(v,x_1)=l+1$ and let $A(x_1)$ be the set of
edges incident with $x_1$ and belonging to $A$. Call $A'$ the
connected normal edge-set resulting from deleting the vertex $x_1$.
By the induction hypothesis, there exits a unique pair $(\psi,
\tilde{B})$ associated to $A'$, then $B=\tilde{B}\cup
\bar{\varphi}(A(x_1))$. Suppose that $B$ is unique for $m=l+1$
and having $n-1$ vertices such that $d(v,x_i)=l+1$. Using the same
argument as in the previous case, it is easy to prove the result
for $n$ vertices.
\end{proof}

Using the same arguments as those given in \cite{mmnr} for locally grid graphs, we prove the
following result.

\begin{lemma}\label{asignarpalabra}
Given a connected normal edge-set $A\subseteq E(H)$, there exits
a unique word $\gamma \in \Gamma$ and an orientation
$(v,v',e,f)$, not necessarily unique, such that
$A=A^{H}(\gamma,v,v',e,f)$.
\end{lemma}

\begin{lemma}\label{versionconexa}
For fixed $n>0$, the number of connected normal edge-sets with
rank $n-1$ and size $n$ is the same for all hexagonal tilings $H$
with $2pq$ vertices and such that $n\leq l_H$.
\end{lemma}

\begin{proof}
By the previous lemma, every connected normal edge-set $A$ is the
instance of a unique word $\gamma(A)$. By Lemma \ref{asignarBaA}
$A$ is isomorphic to $B(\gamma(A))$ then they have the same rank
and size. By Lemma \ref{Bgamma} the instance of a word $\gamma(A)$
is a normal connected edge-set. Hence, the number of connected
normal edge-sets with rank $n-1$ and size $n$ is equal to the
number of distinct instances of words corresponding to edge-sets
of $H^{\infty}$ with rank $n-1$ and size $n$. For every $\gamma$
fixed we can choose $4pq$ different orientations. Let
$\Gamma^{r,s}$ be the set of words of $\Gamma$ such that
$r(B(\gamma)=r$ and $|B(\gamma)|=s$. Then the number of connected
normal edge-sets with rank $n-1$ and size $n$ is:
$$\sum_{\gamma\in \Gamma^{n-1,n}} \frac{4pq}{sym(\gamma)}$$ Since
$sym(\gamma)$ depends only on the word, this quantity does not
depend on the graph.
\end{proof}

A non-connected version of this lemma is proved for hexagonal tilings as  in \cite{mmnr} for
locally grid graphs.

\begin{lemma}\label{versionnoconexa}
For fixed $n>0$ and $t>1$, the number of normal edge-sets with
rank $n-1$, size $n$ and $t$ connected components is the same for
all hexagonal tilings $H$ with $2pq$ vertices and such that $n\leq
l_H$.
\end{lemma}

As a consequence of Lemmas \ref{versionconexa} and
\ref{versionnoconexa}, we have the following result.

\begin{lemma}\label{conexoynoconexo}
Fix $n>0$. The number of edge-sets with rank $n-1$ and size $n$
that do not contain essential cycles is the same for all hexagonal
tilings $H$ with $2pq$ vertices and such that $n\leq l_H$.
\end{lemma}

\begin{corollary}\label{ciclosyTutte}
Let $H$, $H'$ be a pair of hexagonal tilings with $2pq$ vertices.

\noindent $\bullet$ $l_{H}\neq l_{H'}$ implies $T(H;x,y)\neq
T(H';x,y)$.

\noindent $\bullet$ If $l_{H}= l_{H'}$ but $H$ and $H'$ do not
have the same number of shortest essential cycles then
$T(H;x,y)\neq T(H';x,y)$.
\end{corollary}

\begin{proof}
If $l_{H}< l_{H'}$, there are essential cycles of length $l_{H}$
in $H$ but not in $H'$. Since by the preceding lemma the number of
normal edge-sets with rank $l_H-1$ and size $l_H$ is the same in
$H$ and $H'$, the coefficient of $x^{l_H-1}y^{l_H} $ in the
rank-size generating polynomial is greater in $H$ than in $H'$.
In a analogous way the second statement is proved.
\end{proof}

The process we have followed to count normal edge-sets of size
smaller than $l_H+1$, is to show that we can assign a unique set,
$A^{H}(\gamma(A),v,v',e,f)$ to every normal edge-set, $A$. Our
next aim is to prove that for some special cases we can count
normal edge-sets with size greater than $l_H$. This fact is going
to be essential in the next section. Denote by $N(\gamma,H,n,r)$
the number of normal edge-sets in $H$ with rank $r$, size $n$ and
having $B(\gamma)$ as a subgraph. An edge-set $A\subseteq E(H)$
is a \emph{forbidden edge-set for $\gamma$} if it contains an
essential cycle and a subset $B\subseteq A$ isomorphic to
$B(\gamma)$.

\begin{corollary}\label{corolario16}
Given $\gamma \in \Gamma$ such that $B(\gamma)$ contains at least one cycle, then we have that
$N(\gamma,H,n,r)$ is the same for all hexagonal tilings $H$ with $2pq$ vertices, with no forbidden
edge-set for $\gamma$ of size $n$ and such that $n\leq l_H+3$.

\end{corollary}

\begin{proof}
Fix $\gamma \in \Gamma$ such that $B(\gamma)$ contains at least
one cycle. We just have to prove that given $\gamma' \in \Gamma$
such that $B(\gamma)\subseteq B(\gamma')$, $r(B(\gamma'))=l_H+2$
and $|B(\gamma')|\leq l_H+3$, then $A^H(\gamma', v,v',e,f)$ is a
normal edge-set and it has the same rank and size as
$B(\gamma')$. This result is analogous to Lemma \ref{Bgamma} but
taking size smaller than $l_H+4$.

If $B(\gamma)\subseteq B(\gamma')$ then $B(\gamma')$ contains at
least one cycle. If $n\leq l_H$, the proof is equal to the one for
Lemma \ref{Bgamma}. If $n =l_H+1$, we show that $A^H(\gamma',
v,v',e,f)$ has to be an essential cycle (then the rest of the
proof would be analogous to the one followed in Lemma
\ref{Bgamma}). Suppose that $A^H(\gamma', v,v',e,f)$ contains an
essential cycle and one more edge. Since: $$A^H(\gamma,
v,v',e,f)\simeq B(\gamma) \subseteq B(\gamma')\simeq A^H(\gamma',
v,v',e,f)$$ then $A^H(\gamma',v,v',e,f)$ is a forbidden edge-set
for $\gamma$ of size $n$ and we obtain a contradiction. The cases
in which $n=l_H+2$ and $n=l_H+3$ are analogously proved.
\end{proof}

\section{Tutte Uniqueness}

In this section we apply the results obtained in order to prove
the Tutte uniqueness of $H_{k,m,0}$. First we are going to show
that the local hexagonal tiling structure is preserved by the
Tutte polynomial, using the following lemmas.

\begin{lemma}\label{aristaconectividad} \cite{ggmr}
The edge-connectivity of every hexagonal tiling is equal to three.
\end{lemma}

\begin{lemma}\label{conservaparametros} \cite{miernoy} If $G$ is a $2-$connected simple
graph and $H$ is Tutte equivalent to $G$, then $H$ is also simple
and $2-$connected. If $G$ is a $2-$connected simple graph, then
the following parameters are determined by its Tutte polynomial:

(1) The number of vertices and edges.

(2) The edge-connectivity.

(3) The number of cycles of length three, four and five.
\end{lemma}

Next, we are going to prove that there is at least one
coefficient, in which, the Tutte polynomial associated to
$H_{k,m,0}$ and the Tutte polynomial associated to any other
hexagonal tiling not isomorphic to $H_{k,m,0}$, differ. As we have
explained in the first section, the Tutte polynomial associated
to a graph $H$ tell us for every $i$ and $j$ the number of
edges-sets in $H$ with rank $i$ and size $j$. There are two
different types of these edge-sets: normal edge-sets (do not
contain essential cycles) and edge-sets containing at least one
essential cycle. We recall two results from \cite{ggmr}

\begin{lemma}\cite{ggmr}\label{contaresenciales}
Let $G$ be a hexagonal tiling, then the length $l_{G}$ of their
shortest essential cycles and the number of these cycles is:

$$\begin{tabular}{|c|c|l|} \hline $G$ & $l_G$ & \\ \hline
$H_{k,m,r}$ & $\begin{tabular}{c} $2k$ \\ \noalign{\medskip}
$2(m+1)$
 \\
\noalign{\medskip} $2(m+1+r-\lfloor (m+1)/2 \rfloor)$ \\
\noalign{\medskip} $2k$ \\ \end{tabular}$ &
$\begin{tabular}{l}  if $k<m+1$\\
\noalign{\medskip}
if $r < \lfloor (m+1)/2 \rfloor < \lfloor k/2 \rfloor$ \\
\noalign{\medskip} if $ \lfloor (m+1)/2 \rfloor \leq r \leq
\lfloor k/2 \rfloor$
\\
\noalign{\medskip} if $k=m+1$ \\
\end{tabular}$ \\ \hline
$H_{k,m,a}$ & min$(2k,2m+2)$ &  \\ \hline $H_{k,m,b}$ &
min$(2k,2m+2)$ & \\ \hline $H_{k,m,c}$ & min$(k+1,4m+4)$ & \\
\hline $H_{k,m,f}$ & min$(k,4m+8)$ & \\  \hline $H_{k,m,h}$ &
$2k+2$ & \\ \hline $H_{k,m,g}$ & $\begin{tabular}{c}
$2(k-m)-2\lfloor (m+1)/2 \rfloor +3$ \\ \noalign{\medskip} $k+2$
 \\
\noalign{\medskip}  $k+3$ \\ \end{tabular}$ &
$\begin{tabular}{l}  if $k>2m+1$\\
\noalign{\medskip}
if $k\leq 2m+1$ and $k$ odd \\
\noalign{\medskip} if $ k<2m+1$ and $k$ even
\end{tabular}$
\\ \hline \end{tabular}$$

$$\begin{tabular}{|c|c|l|} \hline $G$ & number of essential cycles of length $l_G$ &   \\ \hline $H_{k,m,r}$ &  $\begin{tabular}{c} $m+1$  \\
\noalign{\bigskip} $k\left(
\begin{array}{c} m+1 \\ \lfloor (m+1)/2 \rfloor -r \end{array}
\right)$   \\
\noalign{\bigskip} $k\left(
\begin{array}{c} r+\lfloor (m+1)/2 \rfloor \\ m \end{array}
\right)$
\\
\noalign{\bigskip} $m+1+k\left(
\begin{array}{c} m+1 \\ \lfloor (m+1)/2 \rfloor -r \end{array} \right)$ \\
\end{tabular}$ & $\begin{tabular}{l}  if $k<m+1$\\
\noalign{\bigskip} \noalign{\medskip}
if $r <  \lfloor (m+1)/2 \rfloor < \lfloor k/2 \rfloor$ \\
\noalign{\bigskip} \noalign{\medskip}  if $ \lfloor (m+1)/2
\rfloor \leq r \leq \lfloor k/2 \rfloor$
\\
\noalign{\bigskip}  \noalign{\medskip}  if $k=m+1$ \\
\end{tabular}$
 \\ \hline
$H_{k,m,a}$ &  $\begin{tabular}{c} $m+1$ \\
\noalign{\medskip} $2^{m+1}$   \\
\noalign{\medskip} $2^{m+1}+m+1$  \\
\end{tabular}$ & $\begin{tabular}{l} if $k<m+1$ \\
\noalign{\medskip} if $k>m+1$   \\
\noalign{\medskip} if $k=m+1$  \\
\end{tabular}$
 \\  \hline $H_{k,m,b}$ & $\begin{tabular}{c} $m+1$  \\
\noalign{\bigskip} $2\left(
\begin{array}{c} m+1 \\ (m+1)/2 \end{array} \right)+4\sum_{j=1}^{(m-1)/4}\left(
\begin{array}{c} m+1 \\ (m+1)/2 -2j\end{array} \right) $
\\ \noalign{\bigskip}
$m+1+2\left(
\begin{array}{c} m+1 \\ (m+1)/2 \end{array} \right)+4\sum_{j=1}^{(m-1)/4}\left(
\begin{array}{c} m+1 \\ (m+1)/2 -2j\end{array} \right) $  \\
\end{tabular}$ & $\begin{tabular}{l} if $k<m+1$ \\
\noalign{\bigskip} \noalign{\medskip} if $k>m+1$   \\
\noalign{\bigskip} \noalign{\medskip}  if $k=m+1$  \\
\end{tabular}$
\\  \hline $H_{k,m,c}$ & $\begin{tabular}{c}  $(k/2)\left(
\begin{array}{c} 2m+2 \\ m+1 \end{array} \right) $
\\ \noalign{\bigskip}
$2k $  \\
\end{tabular}$ & $\begin{tabular}{l}  if $4m+4<k+1$   \\
\noalign{\medskip} \noalign{\smallskip}  if $4m+4>k+1$  \\
\end{tabular}$
\\  \hline $H_{k,m,f}$ & $\begin{tabular}{c}  $(k-1)/2\left(
\begin{array}{c} 2m+4 \\ m+2 \end{array} \right) $
\\ \noalign{\bigskip}
$2 $  \\
\end{tabular}$ & $\begin{tabular}{l}  if $4m+8<k$   \\
\noalign{\medskip} \noalign{\smallskip}  if $4m+8>k$  \\
\end{tabular}$
\\  \hline   $H_{k,m,h}$ &  $2^{k+1}$ & \\ \hline
$H_{k,m,g}$ & $\begin{tabular}{c}  2
\\ \noalign{\medskip}
$2(k+2)$ \\ \end{tabular}$ & $\begin{tabular}{l}   if $k\leq 2m+1$ and $k$ odd   \\
\noalign{\medskip}   if $k<2m+1$ and $k$ even  \\
\end{tabular}$ \\ \hline
 \end{tabular}$$
\end{lemma}

\

\begin{lemma}\cite{ggmr}\label{cromatico}
If $G$ is one of the hexagonal tilings, then the chromatic number
of $G$ is given in the following table:

$$\begin{tabular}{|c|c|c|c|c|c|c|c|} \hline $G$ & $H_{k,m,r}$ &
$H_{k,m,a}$ & $H_{k,m,b}$ & $H_{k,m,c}$ &
 $H_{k,m,f}$ & $H_{k,m,g}$ & $H_{k,m,h}$ \\ \hline \noalign{\smallskip} \hline
$\chi (G)$ & 2 & 2 & 2 & 3 & 3 & 3 & 2  \\ \hline
\end{tabular} $$
\end{lemma}

\

\begin{theorem}\label{Tuttetiling}
Let $G$ be a graph Tutte equivalent to a hexagonal tiling $H$,
then $G$ is an hexagonal tiling.
\end{theorem}

\begin{proof}
By Theorem \ref{clasificacion} and Lemma \ref{aristaconectividad}
, we know that $H$ has $2uv$ vertices for some $u$, $v$, $3uv$
edges and edge-connectivity equal to three. By Lemma
\ref{conservaparametros}, $G$ has the same parameters, so $G$ is
cubic, has $2uv$ vertices and $3uv$ edges. Since $H$ is a
hexagonal tiling, $H$ has girth six hence from the previous
lemma, $G$ has girth at least six.

In order to prove that $G$ is a hexagonal tiling, we must show
one more thing: that both graphs have the same number of cycles of
length six (hexagons). Since $G$ and $H$ are Tutte equivalent,
they have the same number of edge-sets with rank five and size
six. We are going to prove that these sets are cycles of length
six.

Let $C\subseteq E(G)$ be an edge-set with rank five and size six. If $y$
is the number of
connected components of $C$, then $5=rg(C)=|V(C)|-y$. Suppose that $C$
is not a cycle, since the
girths of $G$ and $H$ are at least six we have $|V(C_i)|=|A(C_i)|+1$ where $C_i$ is a connected
component of $C$. Hence $|V(C)|=y+6$ and we obtain a contradiction.
\end{proof}

\,

\begin{theorem}\label{uniquetoro}
The graph $H_{k,m,0}$ is Tutte unique for $m\geq 2$ and $k\geq 3$.
\end{theorem}

\begin{proof}
Let $G$ be a graph with $T(G;x,y)=T(H_{k,m,0}; x,y)$.
By Theorem \ref{clasificacion} and Lemma
\ref{cromatico}, $G$ has to be isomorphic to one and only one
of the following graphs:
$H_{k',m',r}$, $H_{k',m',a}$, $H_{k',m',b}$, $H_{k',m',h}$ .
We prove that $G$ is isomorphic to
$H_{k',m',0}$ with $k=k'$ and $m=m'$ by assuming that $G$ is isomorphic
to each one of the
previous graphs and by obtaining a contradiction
in all the cases except in the aforementioned
case. By Lemma \ref{ciclosyTutte}, we obtain $l_{H_{k,m,0}}=l_{G}$
and the number of shortest
essential cycles has to be the same in both graphs.
The process we follow is to compare
$H_{k,m,0}$ with all the graphs given in the previous list.
In the cases in which the minimum
length of essential cycles or the number of cycles of this minimum length
are different, we would
have that the pair of graphs compared are not Tutte equivalent.
For the sake of brevity, we are only going to show
those cases in which these quantities can coincide.

\

{\bf Case 1.} Suppose  $G \cong H_{k',m',r}$ with $r>0$, $k<m+1$
and $k'<m'+1$.

As a result of Lemma \ref{ciclosyTutte}, $k=k'$ and $m=m'$. Our
aim is to prove that the number of edge-sets with rank $2(m+1)-1$
and size $2(m+1)$ is different for each graph. This would lead to a
contradiction since this number is the coefficient of
$x^{2(m+1)-1}y^{2(m+1)}$ of the rank-size generating polynomial.
If $H_{k,m,0}$ has $k\left(
\begin{array}{c} m+1 \\ \lfloor (m+1)/2\rfloor \end{array} \right)+x$
essential cycles of
length $2(m+1)$, then $H_{k,m,r}$ has $x+ \left(
\begin{array}{c} m+1 \\ \lfloor (m+1)/2-r \rfloor \end{array} \right)$ such cycles,
therefore if we show the existence of a bijection between the corresponding edge-sets with rank
$2(m+1)-1$ and size $2(m+1)$ that are not essential cycles, we would have proved the desired
result.

For every $\alpha$ with $0\leq \alpha\leq m-1$, the set of edges that join a vertex $(\alpha,x)$
to $(\alpha+1,x)$ is denoted $E_{\alpha}$. Let $A$ be an edge-set that is not an essential cycle
with rank $2(m+1)-1$ and size $2(m+1)$ in $H_{k,m,0}$. Define $s(A)$ as min$\{ \alpha \in [0,m-1]$
; $ A\cap E_{\alpha}=\emptyset \}$. The minimum always exits for all edge-set contained in
$E(H_{k,m,0})$ or $E(H_{k,m,r})$ that is not an essential cycle. For every $\alpha$ with $0\leq
\alpha\leq m-1$ we define the bijection, $\varphi_r$ between $\{A\subseteq E(H_{k,m,0})|
r(A)=2(m+1)-1, |A|=2(m+1), s(A)=\alpha\}$ and $\{A\subseteq E(H_{k,m,r})| r(A)=2(m+1)-1,
|A|=2(m+1), s(A)=\alpha\}$ as follows:

If $A\subset E(H_{k,m,0})$, $ \varphi_r(A)=\cup\{
\psi(\{(i,j),(i',j')\}); \{(i,j),(i',j')\}\in A\}$
 where
$$\psi(\{(i,j),(i',j')\})= \left\{
\begin{tabular}{lcl} $\{(i,j),(i',j')\}$
& if & $r+1\leq i,i'\leq m$
\\ \noalign{\medskip}
 $\{(i,j-2r),(i',j')\}$  & if &  $i=0$ and $i'=m$
 \\  \noalign{\medskip}
$\{(i,j-2r),(i',j'-2r)\}$ & if & $
i,i'\in [0,r]$ \\
 \end{tabular}\right.$$

\

{\bf Case 2.} Suppose  $G \cong H_{k',m',a}$ with $k<m+1$ and
$k'<m'+1$.

The arguments to prove that both graphs do not have the same Tutte polynomial are analogous to
those developed in the preceding case. We only want to specify two things. The first one is the
number of essential cycles of length $2(m+1)$ in each graph, and the second one, is the bijection
between the edge-sets with rank $2(m+1)-1$ and size $2(m+1)$ that are not essential cycles
contained in each graph.

In $H_{k,m,a}$ we have $k$ exterior edges giving rise to
displaced ladders $i\times j$ with $i=m+1-j$ and $0\leq j\leq
k<m+1$. On the other hand, in $H_{k,m,0}$ each of the $k$ exterior
edges determine a displaced ladder $(m+1)/2 \times (m+1)/2$ for odd $m$
or $(m+2)/2\times m/2$ if $m$ is even. Hence, making use of
the following property:
$$\left( \begin{tabular}{c} $p$ \\
$n$ \\ \end{tabular}\right)<\left( \begin{tabular}{c} $p$ \\
$\lfloor p/2 \rfloor $ \\ \end{tabular}\right) > \left( \begin{tabular}{c} $p$ \\
$l $ \\ \end{tabular}\right) \ {\rm if } \ n< \lfloor p/2 \rfloor
<l$$ the number of essential cycles of length $2(m+1)$ in
$H_{k,m,a}$ is smaller that in $H_{k,m,0}$.

For every $\alpha$ we define $\varphi_r$ between $\{A\subseteq
E(H_{k,m,0})| r(A)=2m+1, |A|=2m+2, s(A)=\alpha\}$ and
$\{A\subseteq E(H_{k,m,a})| r(A)=2m+1, |A|=2m+2, s(A)=\alpha\}$ as
follows:

If $A\subset E(H_{k,m,0})$, $ \varphi_r(A)=\cup\{
\psi(\{(i,j),(i',j')\}); \{(i,j),(i',j')\}\in A\}$
 where
$$\psi(\{(i,j),(i',j')\})= \left\{
\begin{tabular}{lcl} $\{(i,j),(i',j')\}$
& if & $r+1\leq i,i'\leq m$
\\ \noalign{\medskip}
 $\{(i,2k-1-j+3),(i',j')\}$  & if &  $i=0$ and $i'=m$
 \\  \noalign{\medskip}
$\{(i,2k-1-j+3),(i',2k-1-j'+3)\}$ & if & $i,i'\in [0,r]$ \\
 \end{tabular}\right.$$

\

{\bf Case 3.} Suppose  $G \cong H_{k',m',b}$ with $k<m+1$ and
$k'<m'+1$.

As in the previous case, we prove that the number of
essential cycles of length $2(m+1)$ in $H_{k,m,b}$ is smaller than that
in $H_{k,m,0}$.

For every $\alpha$, the bijection $\varphi_r$ between
$\{A\subseteq E(H_{k,m,0})| r(A)=2m+1, |A|=2m+2, s(A)=\alpha\}$
and $\{A\subseteq E(H_{k,m,b})| r(A)=2m+1, |A|=2m+2,
s(A)=\alpha\}$ is defined as follows:

If $A\subset E(H_{k,m,0})$, $ \varphi_r(A)=\cup\{
\psi(\{(i,j),(i',j')\}); \{(i,j),(i',j')\}\in A\}$
 where
$$\psi(\{(i,j),(i',j')\})= \left\{
\begin{tabular}{lcl} $\{(i,j),(i',j')\}$
& if & $r+1\leq i,i'\leq m$
\\ \noalign{\medskip}
 $\{(i,2k-j),(i',j')\}$  & if &  $i=0$ and $i'=m$
 \\  \noalign{\medskip}
$\{(i,2k-j),(i',2k-j')\}$ & if & $i,i'\in [0,r]$ \\
 \end{tabular}\right.$$

 \

{\bf Case 4.} Suppose  $G \cong H_{k',m',a}$ with $k<m+1$ and
$k'>m'+1$.

By Lemma \ref{ciclosyTutte}, $k=m'+1$ and $k'=m+1$. We are going
to prove that there are more edge-sets with rank $2k+2$ and size
$2k+3$ in $H_{k,m,0}$ than in $H_{k',m',a}$. These sets can be
classified into three groups:

{\bf 1.-} Essential cycles of length $2k+3$.

{\bf 2.-} Normal edge-sets (they are edge-sets that do not
contain any essential cycle).

{\bf 3.-} Sets containing an essential cycle of length $2k$ and
three more edges (Figure \ref{normales}).

\begin{figure}[htb]
\begin{center}
\includegraphics[width=70mm]{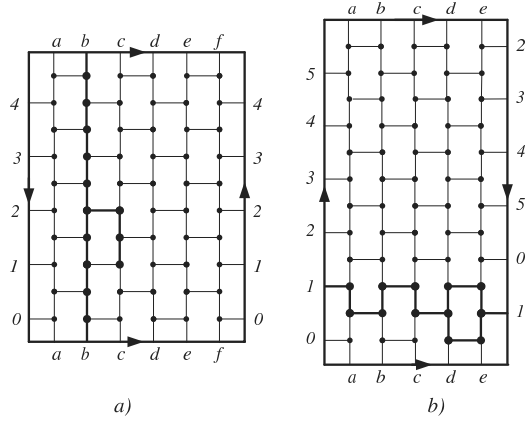} 
\end{center}
\caption{a) Essential cycle of length $2k$ plus three edges in
$H_{k,m,0}$ b) Essential cycle of length $2k=2(m'+1)$ plus three
edges in $H_{k',m',a}$} \label{normales}
\end{figure}

\,

{\bf (1)} Since $H_{k,m,0}$ and $H_{k',m',a}$ are bipartite, they
do not have any essential cycle of odd length.

\,

{\bf (2)} By Corollary \ref{corolario16} we know that $H_{k,m,0}$
and $H_{k',m',a}$ have the same number of normal edge-sets with
rank $2k+2$ and size $2k+3$ that do not contain a cycle of length
six. We are going to prove that the number of normal edge-sets
with rank $2k+2$ and size $2k+3$ containing a cycle of length six
is greater in $H_{k,m,0}$ than in $H_{k',m',a}$.

Again by Corollary \ref{corolario16} , the number of edge-sets
with rank $2k+1$ and size $2k+2$ containing a hexagon is the same
in both graphs, call it $s_{2k+1}$. Add one edge to each of these
sets in order to obtain a set with rank $2k+2$ and size $2k+3$.
This set can be one of the following types depending on which
edge we are adding:

(A) A normal edge set with rank $2k+2$.

(B) A normal edge set containing two non essential cycles and
having rank $2k+1$.

(C) A edge set containing an essential cycle of length $2k$ and a
cycle of length six.

Let $A(H)$, $B(H)$ and $C(H)$ (where $H$ is $H_{k,m,0}$ or
$H_{k',m',a}$) be the number of edge-sets in $H$ that belong to
the groups A, B and C respectively. $A(H)+B(H)+C(H)$ is the
number of possibilities of adding one edge to the edge-sets with
size $2k+2$ containing a hexagon and such that the added edge does
not belong to a fixed hexagon. $$A(H)+B(H)+C(H)=
\frac{s_{2k+1}(3k(m+1)-2k-2)}{2k+3-6}$$ Now, we are going to show
that $(2k-3)B(H)=\sum_{B\in B(H)}(2k+3-\delta (B))$, where
$\delta(B)$ is the number of edges of $B$ which do not belong to
all cycles of length six in $B$. In order to obtain an edge-set
$B\in B(H)$ we add one edge to an edge-set with size $2k+2$ that
contains one hexagon, forming a new contractible cycle. Since the
cycle is contractible, we can consider that the added edge has to
belong to a new hexagon. Hence: $$B(H)= \frac{\sum_{B\in
B(H)}(2k+3-\delta (B))}{2k+3-6}$$ If we apply Corollary
\ref{corolario16} to every $B\in B(H)$, the number of normal
edge-set of size $2k+3-\delta (B)$ is equal in $H_{k,m,0}$ and
$H_{k,m,a}$ then: $$\sum_{B\in B(H_{k,m,0})}(2k+3-\delta
(B))=\sum_{B\in B(H_{k',m',a})}(2k+3-\delta (B))$$
Altogether and by taking in account that $C(H_{k,m,0})=0$
and $C(H_{k',m',a})\neq 0$ (Figure \ref{normales}) we have proved
that $A(H_{k,m,0})> A(H_{k',m',a})$.

\,

{\bf (3)} In $H_{k,m,0}$, every essential cycle of length $2k$
plus three edges has rank $2k+2$ (Figure \ref{normales}a), but in
$H_{k',m',a}$ there are essential cycles in which if we add three
edges we obtain sets with rank $2k+1$ (Figure \ref{normales}b). By
hypothesis, both graphs have the same number of shortest
essential cycles therefore the number of edge-sets in this case
is greater in $H_{k,m,0}$ than in $H_{k',m',a}$.

\

{\bf Cases 5 and 6} If $G\cong H_{k',m',b}$ with $k<m+1$ and
$k'>m'+1$ or $G\cong H_{k',m',h}$ with $k<m+1$ we use the same
arguments as those in the previous case to prove that $H_{k,m,0}$ has
more edge-sets with rank $2k+2$ and size $2k+3$ than
$H_{k',m',b}$ and $H_{k',m',h}$.
\end{proof}

\section{Concluding Remarks}
We have proved that the toroidal hexagonal tiling, $H_{k,m,0}$ is Tutte unique.
Since the results
given in Sections 3 and 4 are proved for all hexagonal tilings,
it seems that all hexagonal
tilings are Tutte unique. To verify this, we would have to do
all cross-checkings between
any two of these graphs as was done in Section 5. This study has been done
for locally grid graphs in \cite{gmr}.

The technique developed in this paper can also be applied to
locally $C_6$ graphs. This would follow from the local orientation existing
in these graphs, based on the minor relationship described in
\cite{ggmr} between locally grid graphs and locally $C_6$ graphs.

\end{document}